\documentstyle[%12pt,
amssymb,amsfonts]{amsart}

\newcommand{\beq}[1]{\begin{equation}\label{#1}}
\newcommand{\enq}[0]{\end{equation}}

\newcommand{\C}[2]{{{#1}\choose{{#2}}}}
\newcommand{\ga}[0]{\alpha }
\newcommand{\gb}[0]{\beta }
\newcommand{\gc}[0]{\gamma }
\newcommand{\gd}[0]{\delta }

\newcommand{\go}[0]{\omega}
\newcommand{\Om}[0]{\Omega}

\newcommand{\gz}[0]{\zeta}
\newcommand{\eps}[0]{\varepsilon }

\newcommand{\gp}[0]{\varphi}

\newcommand{\0}[0]{\emptyset}

\newcommand{\ra}[0]{\rightarrow}
\newcommand{\Ra}[0]{\Rightarrow}

\renewcommand{\lll}[0]{l}

\newcommand{\ttt}[0]{t}

\newcommand{\ww}[0]{{\sf w}}

\newcommand{\E}[0]{{\sf E}}

\newcommand{\A}[0]{{\cal A}}
\newcommand{\B}[0]{{\cal B}}
\newcommand{\cee}[0]{{\cal C}}
\newcommand{\eee}[0]{{\cal E}}
\newcommand{\f}[0]{{\cal F}}

\newcommand{\h}[0]{{\cal H}}

\newcommand{\k}[0]{{\cal K}}

\newcommand{\pee}[0]{{\cal P}}
\newcommand{\r}[0]{{\cal R}}
\newcommand{\s}[0]{{\cal S}}

\newcommand{\Y}[0]{{\cal Y}}
\newcommand{\Z}[0]{{\cal Z}}

\newcommand{\zz}[0]{J}
\newcommand{\TTT}[0]{T}
\newcommand{\RRR}[0]{T}

\newcommand{\mn}[0]{\medskip\noindent}

\newcommand{\sub}[0]{\subseteq}
\newcommand{\sm}[0]{\setminus}
\renewcommand{\dots}[0]{,\ldots,}

\newcommand{\ov}[0]{\overline}

\def\maxr{{\rm maxr \,\,}}
\def\med{{\rm med\,}}

\newcommand{\1}[0]{{\bf 1}}

 \newtheorem{theorem}[subsection]{Theorem}
  \newtheorem{conjecture}[subsection]{Conjecture}

  \newtheorem{lemma}[subsection]{Lemma}
  \newtheorem{corollary}[subsection]{Corollary}
  \newtheorem{definition}[subsection]{Definition}
\theoremstyle{remark}
  \newtheorem{remark}[subsection]{Remark}

\theoremstyle{definition}

\include{psfig}
\begin{document}
\title[Factors in random graphs]
{Factors in random graphs}

\author{Anders Johansson}

\address{Department of Mathematics, Natural and Computer Sciences,
University of G\"avle, SWEDEN}
\email{ajj@@hig.se}

\author{Jeff Kahn}

\address{Department of Mathematics, Rutgers University, Piscataway, NJ 08854 USA}
\email{jkahn@@math.rutgers.edu}
\thanks{J. Kahn  is supported by an NSF Grant.}

\author{Van Vu}
\address{Department of Mathematics, Rutgers University, Piscataway, NJ 08854 USA}
\email{vanvu@@math.rutgers.edu}
\thanks{V. Vu  is supported by an NSF Career Grant.}
\begin{abstract}
Let $H$ be a fixed graph on $v$ vertices.
For an $n$-vertex graph $G$ with $n$ divisible by $v$,
an $H$-{\em factor} of
$G$ is a collection of $n/v$
copies of $H$ whose vertex sets partition $V(G)$.

In this paper we consider the threshold $th_{H} (n)$ of the property that
an Erd\H{o}s-R\'enyi random graph (on $n$ points)  contains an $H$-factor. Our results determine
$th_{H} (n)$ for all strictly balanced $H$.

The method here extends with no difficulty to hypergraphs. As a corollary, we obtain the threshold
for a perfect matching in random
$k$-uniform hypergraph, solving the well-known ``Shamir's problem."

\end{abstract}
\maketitle
\section{Introduction}

%{\bf Thresholds for $H$-factors.}
%
Let $H$ be a fixed graph on $v$ vertices.
For an $n$-vertex graph $G$ with $n$ divisible by $v$,
an $H$-{\em factor} of
$G$ is a collection of $n/v$
copies of $H$ whose vertex sets partition $V(G)$.
(For the purposes of this introduction,
a copy of $H$ in $G$ is
a (not necessarily induced) subgraph of $G$ isomorphic to $H$;
but we will later find it convenient to deal with {\em labeled}
copies.)
{\em We assume throughout this paper
that $v$
divides $n$.}

Let $G(n,p)$ be the Erd\H os-R\'enyi random graph with edge
density $p$.  Recall that
a function $f(n)$ is said to be a {\em threshold} for
an increasing graph property $Q$ if
\begin{equation}\label{thresh}
\begin{array}{ll}
\Pr(\mbox{$G(n,p)$ satisfies $Q$}) \rightarrow 1
&\mbox{if $~p=p(n)=\omega(f(n))$, and}\\
\Pr(\mbox{$G(n,p)$ satisfies $Q$}) \rightarrow 0
&\mbox{if $~p=p(n)=o(f(n))$.}
\end{array}
\end{equation}
Of course if $f(n)$ is a threshold, then so is $cf(n)$ for any
positive constant $c$; nonetheless, following common practice,
we will sometimes say ``{\em the} threshold for $Q$."

Here we are interested in the (increasing) property
that $G(n,p)$ has an
$H$-factor. We denote by $th_H(n) $
a threshold for this property
(where, to be precise, we really mean (\ref{thresh}) holds for
$n\equiv 0 \pmod{v}$).
Determination of the thresholds
$th_H $ is a central problem in the theory of random graphs,
particular cases of which have been
 considered e.g. in
 \cite{Alon-Yuster, Kim, Kriv, Rucinski1}.

In this paper we  come close to a
complete solution to this problem, determining
$th_{H}(n)$ for strictly balanced graphs
(defined below), and getting to within
a factor $n^{o(1)}$ of optimal in general.
Our results generalize without difficulty to hypergraphs.
Here the simplest case---namely, where $H$ consists of
a single hyperedge---settles the much-studied
``Shamir's problem"
on the threshold for perfect
matchings in random hypergraphs of a fixed edge size.
(We will say a little more about this below, but for simplicity
have chosen to give detailed arguments only for strongly
balanced graphs and just comment
at the end of the paper on how these arguments extend.)

 To get some feeling for the magnitude of $th_{H}(n)$, let us first consider lower bounds.
 If $G$ contains an $H$-factor, then each vertex of $G$ is
 covered by some copy of $H$. So, if we denote by
 $th^{[1]}_{H } (n)$ the threshold for the property that in $G(n,p)$
 every vertex is  covered by a copy of $H$, then
\begin{equation}
\label{veryweak} th^{[1]}_{H} (n) \le th_{H} (n).
\end{equation}
For strictly balanced graphs (\ref{veryweak})
will turn out to give the correct value for $th_H$.
But for graphs in which some regions are denser than (or at least
as dense as) the whole, we can run into the problem that, though
the copies of $H$ cover $V(G)$, there is at least one
vertex $x$ of $H$
for which only a few vertices of $G$ play the role of $x$
(in the obvious sense) in these copies.
Formally we have
\begin{equation}
\label{weak} (th^{[1]}_{H} (n) \le) ~\/~
th^{[2]}_{H} (n) \le th_{H} (n),
\end{equation}
where $th_{H} ^{[2]} (n)$ is the threshold for the property:

 \begin{itemize}
\item
every vertex of $G$ is covered by at least one copy of $H$, and

\item for each $x\in V(H)$, there are at least $n/v$ vertices
 $x'$ of $G$ for which some isomorphism of $H$ into $G$ takes
 $x$ to $x'$.
\end{itemize}
(We could replace $n/v$ by $\Omega(n)$ without affecting
this threshold.)

We believe $th^{[2]}_{H}$ tells the whole story:

\begin{conjecture} \label{conj:main}
For every H,
\begin{equation} \label{weak} th^{[2]}_{H} (n) = th_{H} (n). \end{equation}
\end{conjecture}
Though Conjecture \ref{conj:main} may already be a little
optimistic, it does not seem impossible that even a
``stopping time" version is true, {\em viz.}:
if we start with $n$ isolated vertices and add random (uniform)
edges $e_1,\ldots$, then with probability $1-o(1)$ we have
an $H$-factor as soon as we have the property described above in
connection with $th^{[2]}_{H}$.
See \cite{ER} and
\cite{KS,B} for the analogous statements for matchings
and Hamiltonian cycles respectively.

We next want to say something about the behavior of $th^{[1]}_{H}$
and $th^{[2]}_{H}$.  For the issues we are considering
the following notion of
density is natural.
Throughout the paper we use $v(G)$ and $e(G)$ for the numbers
of vertices and edges of a graph $G$.
\begin{definition}
For a graph $H$ on at least two vertices,
$$d(H)= \frac{ e(H)}{ v(H)-1} $$
and
$$ d^{\ast}(H) = \max_{H' \sub H} d(H'). $$
\end{definition}

\begin{definition} A graph $H$ is
{\em strictly balanced} if for any proper
subgraph $H'$ of $H$ with at least two vertices,
$d(H') < d( H).$
\end{definition}
(When the weaker condition $d^{\ast}(H) = d(H)$
holds $H$ is said to be {\em balanced}.)

Strictly balanced graphs are the primary concern of the
present paper and are of considerable
importance for the theory
of random graphs in general; see for example \cite{Boll81,  JLR}.
Basic graphs such as  cliques and cycles are strictly balanced,
as is a typical random graph.

The thresholds $th^{[1]}_{H}$ are known;
see \cite[Theorem 3.22]{JLR}.
In particular, for any $H$ one has
$th^{[1]}_{H} (n) =\Omega (n^{-1/d(H)} (\log n)^{1/m})$, where
$m$ is the number of edges in $H$.
To see why this is natural,
notice that the expected number of copies of $H$ covering a fixed $x\in V(G(n,p))$ is $\Theta(n^{v-1} p^{m})$.
We may then guess that (i) if this expectation is much less
than $\log n$ then the number is zero with probability much
more than $1/n$, and (ii) in this case it is likely that there are
$x$'s for which the number {\em is} zero.

To specify $th^{[2]}_{H}$, we need a little notation.
For $v\in V(H)$ let
$$d^*(v,H)=\max\{d(H'):H'\sub H, v\in V(H')\}$$
(the ``local density" at $v$).
Then clearly $d^*(v,H)\leq d^*(H)$ for all $v$, with equality
for some $v$.
Let $s_v=\min\{e(H'):H'\sub H,v\in V(H'), d(H')=d^*(v,H)\}$
and let $s$ be the maximum of the $s_v$'s.
(Similar notions enter into the determination of $th^{[1]}_{H}$;
again see \cite{JLR},
noting that our uses of $s_v$ and $s$ are not the same
as theirs.)
A proof of the following assertion will appear separately.

\begin{lemma} \label{lemma:lowergeneral}
If $d^*(v,H)=d^*$ $\forall v\in V(H)$ then
$$th_H^{[2]}(n) = n^{-1/d^*}\log^{1/s}n.$$
Otherwise
$$th_H^{[2]}(n) = n^{-1/d^*}.$$
\end{lemma}

In particular we have
\begin{equation}\label{fact:1}
th_H^{[1]}(n)  =th_{H}^{[2]}(n) = n^{-1/ d(H)} (\log n)^{1/m}
\end{equation}
for a strictly balanced $H$ with with $m$ edges
(see \cite{Spencer, Rucinski1}),
and
\begin{equation}\label{fact:2}
th_{H}^{[2]}(n) \ge  n^{-1/ d^{\ast}(H)}
\end{equation}
for a general $H$.

Conjecture \ref{conj:main} says that the values for $th_{H}^{[2]}$
given in Lemma \ref{lemma:lowergeneral}
are also the
values for $th_H$.
Some cases of the conjecture were known earlier:
when $H$ consists of a single edge, it is just the
classic result of Erd\H{o}s and R\'enyi \cite{ER}
giving $\log n/n$ as the
threshold for a perfect matching in $G(n,p)$; and in
\cite{Rucinski1,Alon-Yuster}
(see also \cite[Section 4.2]{JLR})
it is proved whenever $d^*(H)> \gd(H)$ (the minimum degree of $H$).

 The most studied case of the problem is that when $H$ is a
triangle, for which the conjectured value is
$\Theta (n^{-2/3} (\log n)^{1/3} )$.
This problem was apparently first suggested by Ruci\'nski
\cite{Rucinski2}; see \cite[Section 4.3]{JLR} for further
discussion.
Partial results of the form $th_H(n) \le O(n^{-2/3+\eps})$
were obtained by Krivelevich  \cite{Kriv} ($\eps=1/15$)
and Kim \cite{Kim} ($\eps=1/18$).

\section{New results}

The main purpose of this paper is establishing
Conjecture \ref{conj:main} for strictly balanced graphs:

\begin{theorem} \label{theorem:1} Let $H$ be a   strictly balanced
graph with $m$ edges. Then
$$th_H(n)  =  \Theta(n^{-1/ d(H)} (\log n)^{1/m}). $$
\end{theorem}

For general graphs, we have

\begin{theorem} \label{theorem:1-1} For an arbitrary $H$,
$$th_H(n)  =  O(n^{-1/ d^{\ast}(H) +o(1)}). $$
\end{theorem}
This is Conjecture 3.1 of \cite{Alon-Yuster}.
Note that, in view of (\ref{fact:2}),
the bound is sharp up to the $o(1)$ term.

We will actually prove more than Theorems \ref{theorem:1}
and \ref{theorem:1-1}, giving (lower) bounds
on the {\em number} of $H$-factors.
Here, for simplicity,
we confine the discussion to strictly balanced $H$.
From this point through the statement of Theorem
\ref{theorem:3} we
fix  a strictly balanced  graph $H$ on $v$ vertices
and $m$ edges (thus $d(H)= m/(v-1)$).

We use $\f(G)=\f_H(G)$ for the set of $H$-factors in $G$ and
$\Phi(G)= \Phi_H(G)=|\f(G)|$.
For the complete graph
$K_n$ this number is
\begin{equation} \label{equation:F0}
\frac{n ! }{ (n/v) !}  \Big(\frac{1}{|Aut (H)|}
\Big)^{n/v}= n^{\frac{v-1}{v}n}  e^{-O(n) }\end{equation}

\noindent where  $Aut(H)$ is the automorphism group of $H$.
(Note:  in general
the expressions $O(n)$, $o(n)$ can be positive
or negative, so the minus sign is not really necessary here.)
Thus, recalling $\E$ denotes expectation, we have
\begin{equation} \label{equation:expectation}
\E\Phi(G(n,p)) =
n^{\frac{v-1}{v}n} e^{ -O(n) } p^{mn/v} =e^{-O(n)} (n^{v-1} p^m)^{n/v}.
\end{equation}
Strengthening
Theorem \ref{theorem:1}, we will show that, for a suitably large
constant $C$ and
$$p> C n^{-1/d(H)} (\log
n)^{1/m},$$
$\Phi(G(n,p))$ is close to its expectation:

\begin{theorem} \label{theorem:2}
For any $C_1$ there is a $C_2$ such
that for  any
$p > C_2 n^{-1/d(H)} (\log n)^{1/m}$,
$$\Phi(G(n,p)) = e^{-O(n)} (n^{v-1} p^m)^{n/v}$$
with probability at least
$1- n^{-C_1}.$
\end{theorem}
Again note $O(n)$ can be positive or
negative.
The upper bound
%(when $o(n)$ is positive)
follows from
\eqref{equation:expectation} via Markov's inequality. Our task is
to prove the lower bound, for which it will be convenient
to work with the following alternate version.

We say  that an event holds {\em with very high probability (w.v.h.p.)}
if the probability that it fails is $n^{-\omega(1)}$.
We will prove the following equivalent form of Theorem \ref{theorem:2}.
(The equivalence, though presumably well-known, is proved
in an appendix at the end of the paper.)

\begin{theorem} \label{theorem:3}
For $p = \omega (n^{-1/d(H)} (\log n)^{1/m})$,
the  number of $H$-factors in $G(n,p)$ is, with very high
probability, at least $e^{-O(n)} (n^{v-1} p^m)^{n/v}$.
\end{theorem}
As suggested above, the analogous extension of
Theorem \ref{theorem:1-1} also holds.

Finally, we should say something about hypergraphs.
Recall that
a {\em k-uniform hypergraph on vertex set V} is simply
a collection of $k$-subsets, called {\em edges}, of $V$.
Write $\h_k(n,p)$ for the random $k$-uniform hypergraph on vertex
set $[n]$; that is, each $k$-set is an edge with probability $p$,
independent of other choices.

The preceding results extend essentially {\em verbatim} to
$k$-uniform hypergraphs with a fixed $k$.
In particular, as the simplest case (when $H$ is one hyperedge),
we have a resolution of the question, first
studied by
Schmidt and Shamir \cite {SS}
and sometimes referred to as
``Shamir's problem": for fixed $k$ and $n$
ranging over multiples of $k$, what is the threshold for
$\h_k(n,p)$ to contain a perfect matching (meaning, of course, a
collection of edges partitioning the vertex set)?
The problem seems to have first appeared in
\cite{Erdos}, where Erd\H{o}s says that he
heard it from E. Shamir; this perhaps explains the name.

Shamir's problem has been one of the most studied
questions in probabilistic combinatorics over the last 25 years.
The natural conjecture, perhaps first
explicitly proposed in \cite{CFMR},
is that
the threshold is $n^{-k+1}\log n$, reflecting the idea that,
as for graphs, isolated vertices are the primary obstruction
to existence of a perfect matching.
Progress in the direction of this conjecture
and related results
may be found in,
for example, \cite{SS, FJ,Kim,CFMR,Kriv2}.
Again, some discussion of the problem may be found in
\cite[Section 4.3]{JLR}.

We just state the threshold results for hypergraphs, the
relevant definitions and notation extending
without modification to this context.

\begin{theorem} \label{theorem:hyper}
 For a strictly balanced k-uniform
hypergraph $H$ with $m$ edges,
$$th_H(n)  =  \Theta(n^{-1/ d(H)} (\log n)^{1/m}). $$
\end{theorem}

\begin{corollary} \label{cor:hyper}
The threshold for perfect matching in a $k$-uniform random
hypergraph is
 $$ \Theta ( n^{-k+1} \log n ). $$
\end{corollary}
\begin{theorem} \label{theorem:1-1hyper}
For an arbitrary $k$-uniform hypergraph $H$,
$$th_H(n)  =  O(n^{-1/ d^{\ast}(H) +o(1)}). $$
\end{theorem}
Again, the counting versions of these statements also hold
(and are what's actually obtained by following the proof of
Theorem \ref{theorem:2}).

In the  next section  we give  an overview of the proof of
Theorem \ref{theorem:3} and (at the end) an outline of the
rest of the paper.

{\bf Notation and conventions.}
We use asymptotic notation under the
assumption that $n \rightarrow \infty$,
and assume throughout that $n$ is large enough to
support our assertions.
We will often pretend that large numbers are integers,
preferring this common abuse to cluttering the paper
with irrelevant ``floor" and ``ceiling" symbols.

For a graph $G$, $V(G)$ and $E(G)$ denote the vertex and
edge sets of $G$, respectively, and, as earlier, we
use $v(G)$ and $e(G)$ for the cardinalities of these sets.
As earlier, for our fixed $H$ we always take
$v(H)=v$ and $e(H)=m$.

Throughout the paper $V$ is $[n]:=\{1\dots n\}$ and,
as usual, $K_n$ is
the complete graph on this vertex set.
We also use $K_W$ for the complete graph on a general
vertex set $W$.

We use $\log$ for natural logarithm,
$\1_{\eee}$ for the indicator of event $\eee$,
and $\Pr$ and $\E$
for probability and expectation.

As noted earlier, we will henceforth use {\em copy of H in G}
to mean an injection $\varphi:V(H)\ra V(G)$ that takes edges
to edges
(i.e. $\varphi(E(H))\sub E(G)$ with the natural interpretation
of $\varphi(E(H))$).
This avoids irrelevant issues involving $|{\rm Aut}(H)|$
and will later make our lives somewhat easier in other ways
as well.
Note that since this change multiplies $\Phi(G)$ by
$|Aut(H)|^{n/v} = e^{O(n)}$, it does not affect the statements of
Theorems \ref{theorem:2} and \ref{theorem:3}.
Note also that we will still sometimes use (e.g.) $K$ for a copy
of $H$, in which case we think of a labeled copy in the usual sense.

We use $\h(G)$ for the set of copies of $H$ in $G$,
$\h(x,G) =\{K\in \h(G):x\in V(K)\},$
(where $x\in V(G)$), and $D(x,G)=|\h(x,G)|$.

\section{Outline}\label{Outline}

In this section we
sketch the proof of
Theorem \ref{theorem:3}, or, more precisely, give the proof
modulo various assertions (and a few definitions)
that will be covered in later sections.
For the moment we actually work with
the model $G(n,M)$---that is the graph chosen uniformly
from $M$-edge graphs on $V$---though
in proving the assertions made here it will turn out to be
more convenient to return to $G(n,p)$.
The relevant connection between the two models is
given by Lemma \ref{lemma:twomodels}, which in particular
implies that
Theorem \ref{theorem:3} is equivalent to
\begin{theorem} \label{theorem:3'}
For $p=p(n) = \omega (n^{-1/d(H)} (\log n)^{1/m})$
and $M=M(n)=\C{n}{2}p$,
$$\Pr(\Phi(G(n,M) \ge (n^{v-1} p^m)^{n/v}e^{-O(n)} ) \ge  1-n^{-\go(1)}.$$
\end{theorem}
For the proof of this, let $M=M(n)$ be as in the statement
and $\TTT=\C{n}{2}-M$.
Let $e_1\dots e_{\C{n}{2}}$ be a random (uniform) ordering of $E(K_n)$
and set $G_i = K_n-\{e_1\dots e_i\}$ (so $G_0=K_n$)
and $\f_i=\f(G_i)$
(recall this is the set of $H$-factors in $G_{i}$).
Let $\xi_i$
be the fraction of members of $\f_{i-1}$
containing $e_i$ (where we think of an $H$-factor as a subgraph
of $G$ in the obvious way).
Then (for any $t$)

$$| \f_{t}| = |\f_{0}| \frac{|\f_{1}| }{| \f_{0}|} \cdots \frac{ |\f_{t}| }{| \f_{t-1}|} = |\f_{0}|
(1-\xi_{1}) \cdots (1-\xi_{t}), $$
or
\begin{equation}\label{mgt1}
\log |\f_t| =\log|\f_0| +\sum_{i=1}^t \log (1-\xi_i).
\end{equation}
Now
\begin{equation}\label{mg0}
\log|\f_0| =\log
%\left[
\frac{n!}{(n/v)!|{ Aut}(H)|^{n/v}}
%\right]
= \frac{v-1}{v}n\log n -O(n).
\end{equation}
We also have
\begin{equation}\label{gci}
\E \xi_i = \frac{mn/v}{\C{n}{2}-i+1} =:\gc_i,
\end{equation}
since in fact
\begin{equation}\label{infact}
\E [\xi_i|e_1\dots e_{i-1}] =\gc_i
\end{equation}
for any choice of $e_1\dots e_{i-1}$.
Thus
\begin{equation}\label{Em}
\sum_{i=1}^t\E \xi_i = \sum_{i=1}^t \gc_i=
\frac{mn}{v} \log \frac{\C{n}{2}}{\C{n}{2}-t} +o(1),
\end{equation}
provided $\C{n}{2}-t>\omega(n)$.
Let $\A_t $ ($=\A_t(n)$) be the event
$$
 \{\log |\f_t| >\log|\f_0| - \sum_{i=1}^t\gc_i -O(n)\}.
$$
Our basic goal is to show that failure of $\A_t$ is unlikely. Precisely, we want to show
\begin{equation}\label{Main1}
\mbox{for $t\leq \TTT$, $\,\,\,\Pr(\bar{\A}_t) =  n^{-\go(1)}$.}
\end{equation}

This implies  Theorem~\ref{theorem:3'},
since if $\A_T$ holds then (\ref{mg0}) and
(\ref{Em}) (the latter with $t=T$) give
$$
\log \Phi (G_{T}) =\log |\f_T|  > \frac{v-1}{v}n\log n +\frac{mn}{v}\log p  -O(n) .
$$

The proof of (\ref{Main1})
uses the method of martingales with bounded differences (Azuma's inequality).
Here it is  natural to consider the martingale
$$X_t =\sum_{i=1}^t (\xi_i -\gc_i)$$
(it is a martingale by (\ref{infact})), with associated
difference sequence
$$
Z_i = \xi_i -\gc_i.
$$

Establishing concentration of $X_t$ depends on maintaining some
control over the $|Z_i|$'s,
for which purpose we will keep track of
two sequences of auxiliary
events, $\B_i$ and $\r_i$ ($1\le i \le \TTT-1$).
Informally, $\B_i$ says that no copy of $H$ is in much more than its
proper share of the members of $\f_i$, while $\r_i$ includes
regularity properties of $G_i$, together with some
limits on the numbers of copies of fragments of $H$ in $G_i$.
The actual definitions of $\B_i$ and $\r_i$ are given in
Section~\ref{BandR}.

For $i\leq T $
it will follow easily from $\B_{i-1}$ and $\r_{i-1}$
(see Lemma  \ref{proofofxibd}) that
\begin{equation}\label{xibd}
\xi_i = o(\log^{-1} n).
\end{equation}
(The actual bound on $\xi_i$
implied by $\B_{i-1}$ and $\r_{i-1}$
is better when $i$ is not close to $\TTT$, but for our
purposes the difference is unimportant.)

The bound (\ref{xibd}), while (more than) sufficient
for the concentration we need, can occasionally fail, since the auxiliary properties $\B_{i-1}$ and $\r_{i-1}$ may fail.
To handle this problem we slightly modify the preceding
$X$'s and $Z$'s by setting
\begin{equation}\label{Zi}
Z_i =\left\{\begin{array}{ll}
\xi_i -\gc_i & \mbox{if $\B_j$ and $\r_j$ hold for all $j<i$}\\
0&\mbox{otherwise}
\end{array}\right.
\end{equation}
(and
$X_t = \sum_{i=1}^tZ_i$).

As shown in Section \ref{Mart},
%%A
a martingale analysis along the lines of
Azuma's Inequality then gives (for example)
\begin{equation}\label{conc}
\Pr( |X_t| > n  ) < n^{-\omega(1)}.
\end{equation}
Notice that if we do have $\B_i$ and $\r_i$ for $i<t\leq T$
(so that $X_t=\sum_{i=1}^t(\xi_i-\gc_i)$)
and
$|X_t| \leq O(n)$
(it will actually be smaller)
then we have $\A_t$,
since:
$$
\sum_{i=1}^t\xi_i < \sum_{i=1}^t\gc_i + O(n) < O(n\log n),
$$
whence (noting we have
(\ref{xibd}) for $i\leq t$)
$$
\sum_{i=1}^t\xi_i^2 < o(\log^{-1} n)\sum_{i=1}^t \xi_i < o(n),
$$
so that (using (\ref{mgt1}))
$$
\log |\f_t|
> \log|\f_0| -\sum_{i=1}^t(\xi_i +\xi_i^2 ) \\
> \log|\f_0| -\sum_{i=1}^t\gc_i  - O(n).
$$

Thus the first failure (if any) of an $\A_t$ (with $t\leq \TTT$)
must occur either because $X_t$ is too large
or because one of the properties $\B$, $\r$ fails even earlier
(that is, $\B_i$ or $\r_i$ fails for some $i<t$).
Formally we may write
$$
\Pr(\bar{\A}_t) < \sum_{i<t}\Pr(\bar{\r}_i) +
\sum_{i\leq t}\Pr(\wedge_{j<i}(\B_j\r_j)\wedge \bar{\A}_i) +
\sum_{i< t}\Pr(\A_i\r_i\bar{\B}_i).
$$
Here the second sum is, in view of the preceding discussion
and (\ref{conc}), at most
$n^{-\omega(1)}$,
and we will show, for $i\leq \TTT$,
\begin{equation}\label{Ri}
\Pr(\bar{\r}_i) < n^{-\go(1)}
\end{equation}
and
\begin{equation}\label{Bi}
\Pr(\A_i\r_i\bar{\B}_i) < n^{-\go(1)}.
\end{equation}
Combining these three bounds gives
(\ref{Main1}) and Theorem~\ref{theorem:3'}.\qed

The rest of the paper is organized as follows. The next few sections
(\ref{Models}-\ref{Ent}) fill in some general results, and, as
already mentioned, Section \ref{Mart} gives the calculation needed
for (\ref{conc}).  In Section \ref{BandR} we finally define the
properties $\B_i$ and $\r_i$, and slightly restate our main
inequalities, (\ref{Ri}) and (\ref{Bi}). Proofs of these---the heart
of our argument---are given in Sections \ref{Reg}-\ref{PLC}.
Extensions to general graphs and hypergraphs are discussed in
Section \ref{extension}. Finally,  Section \ref{Appendix} gives the
equivalence of Theorems \ref{theorem:2} and \ref{theorem:3}.
% and Section \ref{appendixB} fills in some proofs omitted from
%Section \ref{Conc}.

\section{Models of random graphs} \label{Models}

As noted earlier, we go back and forth between the models
$G(n,p)$ and $G(n,M)$.

We recall the definitions, which were already used above:
$G(n,p)$ is the random graph on vertex set $[n]$
gotten by taking each pair of
vertices to be an edge with probability $p$, independently
of other choices; and
$G(n,M)$ is the random graph on $[n]$ gotten by choosing $M$
edges uniformly from
$E(K_n)$.
% (equivalently, by deleting a uniformly chosen
%$(T:= {n \choose 2}-M)$-subset from $E(K_n)$).)

The following observation allows us to switch
between the two models.
(Note that $ { n \choose 2} p$ is the
expected number of edges in $G(n,p)$.)

\begin{lemma} \label{lemma:twomodels}
Let $n^{\Omega(1)}= M \le {n \choose 2}$ be an integer,
and $p= M/{n \choose 2}$.

{\rm (a)} If an event $\eee$
holds with probability at least $1-\eps$ in $G(n,p)$,
then it  holds with probability at least $1- O(n\eps)$ in $G(n,M)$.

{\rm (b)} If  an event
$\eee$ holds with
probability at least $1-\eps$ in $G(n,M')$, for all
$M/2 \le M' \le 2M$,
then it holds with probability at least
$1- \eps-n^{-\omega (1)}$ in $G(n,p)$.
\end{lemma}

Notice that the error terms in the conclusions are not the
same as those in the hypotheses.  For our purposes the extra
factors will not matter, since
we will always be dealing with failure probabilities of order
$n^{-\omega (1)}$.

{\em Proof of Lemma} \ref{lemma:twomodels}.
The number of edges
in $G(n,p)$ has the binomial distribution ${\rm Bin}( {n \choose 2}, p)$.
Claim (a) follows from the fact that this number is
equal to $M$ with probability
$\Omega (1/n)$, and claim (b) from the fact that it
is between $M/2$ and $2M$ with probability $1- n^{-\go(1)}$.

Both these facts are simple properties of the distribution
${\rm Bin}({n \choose 2}, p)$ for $p$ as in the lemma.

\section{Concentration}\label{Conc}

Let $f=f(t_1 \dots t_n)$ be a polynomial of degree $d$ with real
coefficients. We say that $f$ is {\it normal} if its coefficients
are positive numbers and the maximum coefficient is 1.
We will be interested in the situation where the $t_i$'s
are independent
Bernoulli random variables
(which will be i.i.d. in our applications, but need not be
for the results stated here).
Then $f$ is itself a r.v., and we will want to say that,
under appropriate conditions, it is
likely to be close to its mean.

All polynomials needed for the proof of Theorem~\ref{theorem:3}
are multilinear, that is, of the form
$f(t) =\sum \ga_{_U}t_{_U}$,
where $U$ ranges over subsets of $[n]$ and
$t_{_U}:=\prod_{u\in U}t_u$.
So for simplicity
we confine the present discussion to such polynomials,
though the results hold in greater generality.

For a (multilinear) polynomial $f$ as above and $L\sub [n]$,
the partial derivative (of order $|L|$) with respect to
the
variables indexed by $L$ is
$
\sum_{U\supseteq L}\ga_{_U}t_{U\sm L},
$
and its expectation, denoted $\E_L$ or $\E_Lf$, is
$
\sum\{\ga_U\prod_{i\in U\sm L}p_i:U\supseteq L\},
$
where $t_i\sim {\rm Ber}(p_i)$.
In particular, when $p_i=p ~\forall i$ we have
$
\sum\{\ga_Up^{|U\sm L|}:U\supseteq L\}.
$
Set $\E_j f=\max_{|L|=j}\E_L f$.

%%%%%%%%%%%%%%%%%%%%%%%%%%%%%%%%%%%%%%%%%%%%%%%??
It may be helpful to observe that when, as here, $f$ is
multilinear, we don't really need the polynomial language:
We have a weight function $\ga$ on subsets of $[n]$
and are choosing a random subset, say $T$, of $[n]$,
with membership in $T$ determined by independent coin flips.
Then $t$ is the indicator of $T$;
$f(t)$ is the (total) weight of the surviving sets
(i.e. those contained in $T$); and $\E_Lf$ is the expected
weight of the surviving sets containing $L$, given that
$L$ itself survives.
%%%%%%%%%%%%%%%%%%%%%%%%%%%%%%%%%%%%%%%%%%%%%%%??

We begin with a special case  of the main  result of  \cite{KV}.

\begin{theorem}  \label{theorem:KV}
The following holds for any fixed positive integer $d$ and positive
constant $\eps$. Let $f$ be a multilinear, homogeneous, normal
polynomial of degree $d$ such that $\E f \ge n^{\eps} \max_{1 \le j
\le d} \E_j f$. Then
$$\Pr ( |f-\E f| > \eps \E f) = n^{-\omega (1) }. $$
\end{theorem}
Note that (here and below) $d$ and $\eps$ are fixed,
and expressions such as $\go(1)$ are interpreted as $n\ra \infty$.
\begin{remark}\label{almostnormal}
Though we have stated the results in this section for
normal polynomials, they clearly remain true if we
instead assume some fixed bound on coefficients---in this
case we say a polynomial is $O(1)$-{\em normal}---and
we will sometimes use them in this form.
\end{remark}

When $f$ is normal, $\E_d f=1$; so in order to use
Theorem \ref{theorem:KV}
we must have $\E f = \Omega (n^{\eps})$ (if one follows \cite{KV}
more closely, then $\E f =\omega ((\log n)^{2d})$ is sufficient).
Thus, the theorem is not applicable when, for
example, $\E  f$ is on the order of $\log n$, a range
that will turn out to be crucial for us.
In such ranges
the following theorem, which is essentially (the
multilinear case of)
Theorem 1.2 of \cite{Vu},
is useful.

\begin{theorem} \label{theorem:V}
The following holds for any fixed positive integer d and positive
constant $\eps$. Let $f$ be a multilinear,  normal, homogenous
polynomial of degree $d$ such that $ \E f=\omega (\log n)$  and $
\max_{1 \le j \le d-1} \E_j  f \le n^{-\eps}$. Then
$$\Pr ( |f-\E f| > \eps \E f) = n^{-\omega (1) }. $$
\end{theorem}
%%%%%%%%%%%%%%%%%%%%%%%%%%%%%%%%%%%%%%%%%%%
(Note:  Theorem 1.2 of \cite{Vu} asks that $\E f$ be $o(n)$, but
this requirement is superfluous:  since $f$ is normal, $\E f \le
n^d$; so we can add, say, $n^{d+1}$ dummy variables that don't
appear in $f$ to guarantee $\E f =o(n)$,
and then apply the theorem
as stated. Also, to be very picky, the statement in
\cite{Vu} is in finitary form---for any given $\gb$ the probability is
less than $n^{-\gb}$ if $n/Q > \E f > Q\log n$ for a large enough
$Q=Q_{\gb}$---but this clearly implies the present version.)
%%%%%%%%%%%%%%%%%%%%%%%%%%%%%%%%%%%%%%%%%%%%

The key difference between Theorem \ref{theorem:V}  and Theorem
\ref{theorem:KV} is that here we only need to consider partial
derivatives of order up to $d-1$. While this may seem to be a minor
point, the proof of Theorem \ref{theorem:V} is considerably more
involved than that of Theorem \ref{theorem:KV}.

Our basic concentration  statement is
the following combination of Theorems
\ref{theorem:KV} and \ref{theorem:V}.

\begin{theorem} \label{theorem:combination}
The following holds for any fixed positive integer $d$ and positive
constant $\eps$. Let $f$ be a multilinear, homogeneous, normal
polynomial of degree $d$ such that $ \E f=\omega (\log n)$  and $
\max_{1 \le j \le d-1} \E_j  f \le n^{-\eps}\E f$. Then
$$\Pr ( |f-\E f| > \eps \E f) = n^{-\omega (1) }. $$
\end{theorem}

\begin{corollary} \label{Vucor}
The following holds for any fixed positive integer $d$ and positive
constant $\eps$. Let $f$ be a multilinear, homogeneous, normal
polynomial of degree $d$ with $ \E f\leq A$, where $A$ ($=A(n)$)
satisfies
$$A \geq \omega(\log n) + n^{\eps} \max_{0<j<d}\E_j f ,$$
then
$$\Pr ( f>(1+\eps)A ) \le n^{-\omega(1) }. $$
\end{corollary}
{\em Proof (sketch).}
If $\E f\geq A/2$ then this is immediate from
Theorem \ref{theorem:combination}.
Otherwise we can augment $f$ to some $g\geq f$ for which
we still have $\E h \leq A$ and $g$ satisfies the
hypotheses of Theorem \ref{theorem:combination}.
(For instance, we can do this with $g=g(t,s) =f(t) + h(s)$, where
$h(s) = h(s_1\dots s_n)=\gc \sum\{s_U:U\sub [n], |U|=d\}$
with the $s_i$'s all ${\rm Ber}(p)$
(and different from the $t_i$'s),
and $\gc$ chosen so that $\E h = A/2$.)\qed

We also need an inhomogeneous version of
Theorem~\ref{theorem:combination}.
Write $\E_L' = E_L'f$
for the expectation of the {\em nonconstant} part of the
partial derivative of $f$ with respect to $L$.
Of course for $f$ homogeneous of degree $d$ and $0<|L|<d$
as in Theorem~\ref{theorem:combination}, we have $E_L'f=E_Lf$.

\begin{theorem} \label{inhomog}
The following holds for any fixed positive integer $d$ and positive
constant $\eps$. Let $f$ be a multilinear, normal polynomial of
degree at most $d$, with $ \E f = \omega(\log n)$ and $
\max_{L\neq\0}\E_L' f \le n^{-\eps}\E  f$. Then
$$\Pr ( |f-\E f| > \eps \E f) \le n^{-\omega(1)}. $$
\end{theorem}
{\em Proof (sketch).}
We apply Theorem \ref{theorem:combination}
to each homogeneous part of $f$.
Say the $k$th homogeneous part is $f_k$.
If $\E f_k > \E f/2$ then Theorem \ref{theorem:combination}
gives
$\Pr(|f_k-\E f_k| > (\eps/d)\E f) = n^{-\go(1)}$
(since $\E f_k =\Theta(\E f$)).
Otherwise we can, for example, introduce dummy variables
$z_{ij}$ for $1\leq i \leq \lceil \E f\rceil$ and $j \in [k]$,
each equal to 1 with probability 1,
and let
$g_k = \sum_i \prod_j z_{ij}$ and
$h_k = f_k +g_k$.
Theorem \ref{theorem:combination} (applied to $h_k$, for
which we have $\E h_k =\Theta(\E f$)) then gives
$$\Pr(|f_k - \E f_k|> \frac{\eps}{d}\E f) =
\Pr(|h_k - \E h_k| >\frac{\eps}{d}\E f)\\
= n^{-\go(1)}.
$$
(The number of variables used for $h_k$ is not exactly $n$, but
cannot be big enough to make any difference.) Consequently, we have
$$\Pr ( |f-\E f| > \eps \E f) \leq \sum_{k=1}^d \Pr(|f_k - \E f_k|> \frac{\eps}{d}\E f) =n^{-\omega(1)}.$$\qed

The corresponding generalization of
Corollary~\ref{Vucor} is
\begin{corollary} \label{Vucor2}
The following holds for any fixed positive integer $d$ and positive
constant $\eps$. Let $f$ be a multilinear, normal polynomial of
degree at most $d$ with $ \E f\leq A$, where $A=A(n)$ satisfies
$$A \geq \omega(\log n) + n^{\eps} \max_{L\neq \0}\E_L' f .$$
Then
$$\Pr ( f>(1+\eps)A ) = n^{-\omega(1) }. $$
\end{corollary}

Finally, we will sometimes need to know something when $\E f$ is
smaller.  Here we cannot expect that $f$ is
(w.v.h.p.) close to its mean, but will be able to
get by with the following weaker statement, which is
(in slightly different language)
Corollary 4.9 of
\cite{Vu2}.  Note that in this case our hypothesis
includes $\E_{\0}f$ ($=\E f$).
\begin{theorem}\label{lastVu}
The following holds for any fixed positive integer $d$ and positive
constant $\eps$. Let $f$ be a multilinear, normal polynomial of
degree at most $d$ with $ \max_L\E_L' f  \le n^{-\eps}$. Then for
any $\gb(n) =\omega(1)$,
$$\Pr ( f > \gb(n)) = n^{-\omega(1)}. $$
\end{theorem}

For a random graph $G(n,p)$, we use $t_{e}$ for
the Bernoulli variable representing the appearance of the edge $e$;
thus $t_{e}$ is 1 with probability $p$ and $0$ with
probability $1-p$. Recall that for $S\sub E(K_n)$,
$t_{S} := \prod_{e \in S}  t_{e}. $

To get a feel for where this is headed, let us briefly
discuss a typical use of the preceding theorems.
Consider a  strictly balanced graph $H$ with $v$ vertices and $m$
edges, and suppose
we are interested in the number of copies of $H$ in $G(n,p)$
containing a fixed vertex $x_0$. This number is naturally expressed
as a
multilinear, normal, homogeneous polynomial of degree $m$ in the
$t_e$'s:  $f = \sum_{U} t_U $, where $U$ runs over
edge sets of copies of $H$ in $K_n$ containing $x_0$. Clearly
$\E f= \Theta (n^{v-1} p^m)$.

Now consider a non-empty subset $L$ of $E(K_n)$.
The partial derivative $\E_Lf$
is $\sum_{U\supseteq L} t_{U \backslash L}$.
Let $H'$ be the subgraph of $K_n$ consisting of $L$
and those vertices incident with edges of $L$,
and let $v'$ and $m'$ ($=|L|$) be the numbers of
vertices and edges of $H'$.
Then $\E_Lf$ is $O (n^{v-v'} p^{m-m'})$
if $L$ covers $x_0$ and
$O (n^{v-v'-1} p^{m-m'})$ if it does not.
In either
case
$$\E f/ \E_L f = \Omega (n^{v'-1 } p^{m'}). $$

Since $H$ is strictly balanced, $(v'-1)/m' > (v-1)/m $. This implies
that if $\E f= \Omega (1)$ (that is,
if $p=\Omega(n^{-(v-1)/m}$) then $\E f/ \E_L f \ge
n^{\Omega(1)}$, a hypothesis for most of the theorems
of this section.
(This discussion
applies without modification to hypergraphs, since it makes no
use of the fact that edge are of size two.)

\section{Entropy}\label{Ent}

We use $H(X)=H_e(X)$ for the base $e$
entropy of a discrete r.v. $X$;
that is,
$$H(X) =\sum_x p(x)\log \frac{1}{p(x)},
$$
where $p(x) = \Pr(X=x)$.
(For entropy basics see e.g. \cite{Csiszar-Korner} or \cite{McE}.)

Given a graph $G$ and $y\in V(G)$, we use
$X(y,G)$ for the copy of
$H$ containing $y$ in a uniformly chosen
$H$-factor of $G$, and
$h(y,G)$ for $H(X(y,G))$.
(We will not need to worry about $G$'s without $H$-factors.)
The next lemma is a special case of a fundamental observation
of Shearer \cite{CFGS}.

\begin{lemma}\label{Shearer}
For any G,
$$
\log \Phi(G)\leq \frac{1}{v}\sum_{y\in V(G)}h(y,G).
$$
\end{lemma}
%%%%%%%%%%%%%%%%%%%%%%%%%%%%%%%%%%%%%%%%%%%%%%%%%%
(To be precise, the statement we are using is more general
than what's usually referred to as ``Shearer's Lemma,"
but is what Shearer's proof actually gives; it may
be stated as follows.
Suppose $Y=(Y_i:i\in I)$ is a random vector and
$\s $ a collection of subsets of $I$
(repeats allowed) such that each $i\in I$ belongs
to at least $t$ members of $\s$.
Then $H(Y)\leq t^{-1}\sum_{S\in \s}H(Y_S)$, where
$Y_S $ is the random vector $(Y_i:i\in S)$.
To get Lemma \ref{Shearer} from this, let
$Y$ be the indicator of the random $H$-factor
(so $I$ is the set of copies of $H$ in $K_n$)
and $\s = (S_v:v\in V)$, where $S_v$ is the set
of copies of $H$ (in $K_n$) containing $v$.)
%%%%%%%%%%%%%%%%%%%%%%%%%%%%%%%%%%%%%%%%%%%%%%%%%%

For the next lemma,
$S$ is a finite set, $\ww:S\ra \Re^+$,
and
$X$ is the random variable
taking values in $S$ according to
$$\Pr(X=x) = \ww(x)/\ww(S)$$
(where for $A\sub S$, $\ww(A) = \sum_{x\in A}\ww(x)$).
For perspective recall that for any r.v. $X$ taking values in $S$,
one has $H(X)\leq \log|S|$ (with equality iff $X$ is uniform
from $S$).
\begin{lemma}\label{entlemma}
If $H(X) > \log |S| -O(1)$,
then there are $a,b\in {\rm range}(\ww)$ with
\begin{equation}\label{ab}
a\leq b < O(a)
\end{equation}
such that for $J=\ww^{-1}[a,b]$ we have
$$
|J| = \Omega (|S|)
$$
and
$$
\ww(J) > .7\ww(S).
$$
\end{lemma}
%{\em [Remark.
%If we assume $H(X) >\log |S| -\eps$ then we can
%(I think) strengthen the conclusions to
%$a\leq b<(1+f(\eps))a$,
%$|J|> (1-f(\eps))|S|$ and
%$\ww(J)> (1-f(\eps))\ww(S)$,
%where $f(\eps)\ra 0$ as $\eps \ra 0$.]}

\mn
{\em Proof.}
Let
\begin{equation}\label{H(X)}
H(X)= \log |S|-K
\end{equation}
and define $C$ by
$\log C = 4(K+\log 3)$.
With $\ov{\ww}=\ww(S)/|S|$, let $a=\ov{\ww}/C$,
$b = C\ov{\ww}$,
$L=\ww^{-1}([0,a))$,
$U=\ww^{-1}((b,\infty])$,
and $J= S\sm (L\cup U)$.
We have
\begin{equation}\label{H(X)'}
H(X) \leq \log 3 + \frac{\ww(L)}{\ww(S)}\log |L| +
\frac{\ww(J)}{\ww(S)}\log |J|
+\frac{\ww(U)}{\ww(S)}\log |U|.
\end{equation}
%{\em [Remark.  This proof is surely not optimal.
%For one thing we should probably be writing $\log|S|-H(X)$
%exactly as the divergence of uniform dist. and the
%dist. given by $\ww$.  But for now just getting some version down.]}

Then we have a few observations.
First, $|U|< |S|/C$ implies that the r.h.s. of (\ref{H(X)'})
is less than
$$
\log 3 +\log |S| -\frac{\ww(U)}{\ww(S)}\log C,
$$
which with (\ref{H(X)}) implies
\begin{equation}\label{wU}
\ww(U) < \frac{K+\log 3}{\log C}\ww(S) = \ww(S)/4.
\end{equation}
Of course this also implies $|U|< |S|/4$.

Second, combining (\ref{wU}) with the trivial $\ww(L)< \ww(S)/C$,
we have (say)
$\ww(J)>.7\ww(S)$.
But then (third) since the r.h.s. of (\ref{H(X)'}) is
at most
$$\log 3 +\log |S| +\frac{\ww(J)}{\ww(S)}\log \frac{|J|}{|S|}
<\log 3 +\log |S| +.7\log \frac{|J|}{|S|},
$$
we have
$$
|J|\geq \exp[ -(.7)^{-1}(K+\log 3)]|S| ~(=\Omega(|S|)).
$$

\section{Martingale}\label{Mart}

Here we give the proof of   (\ref{conc}). As in Section 3, let
$$X_{t} = Z_{1} + \cdots +Z_{t}, $$
where we use the modified $Z_i$'s of (\ref{Zi}). We first prove that

\begin{equation} \label{eqn:concentration} \Pr (X_t \ge n) <
n^{-\omega (1)}. \end{equation}
Notice that  $Z_{i}$ is a function of the random sequence $e_{1}
\dots e_{i}$, with $\E(Z_{i}| e_{1} \dots e_{i-1})=0$ for any $e_{1}
\dots e_{i-1}$ and (relaxing (\ref{xibd}) a little)
$|Z_{i}| <\eps :=\log^{-1} n$.
%For convenience, in the rest of this section we set $\eps= \log^{-1}n$.
By Markov's inequality, we have, for any positive $h$,

\begin{equation} \label{eqn:concentration0} \Pr(X_t \ge n) = \Pr(e^{h(Z_1+ \cdots +Z_t)} \ge e^{hn}) \le
\E(e^{h(Z_1+ \cdots +Z_t)})e ^{-hn}. \end{equation}

Next we  bound $\E(e^{h(Z_1+ \cdots +Z_t)})$.
Since $Z_i= \xi_i- \gamma_i$, $\E(\xi_i|e_1 \dots e_{i-1})
=\gamma_i$ and $0\leq \xi_i \le \eps$, we have (using
convexity of $e^x$),
$$\E( e^{hZ_i}| e_1 \dots e_{i-1}) \le e^{-h\gamma_i}
((1- \frac{\gamma_i}{\eps}) + \frac{\gamma_i}{\eps} e^{h \eps} ). $$

A simple Taylor series calculation shows that the right hand side is
at most $e^{h^2 \eps \gamma_i}$, for any $0\le h \le 1$. Thus, for
such $h$
%
%\begin{equation}\label{eqn:concentration1}
$$ \E( e^{hZ_i}| e_1 \dots e_{i-1}) \le e^{h^2 \eps \gamma_i},$$
%\end{equation}

and induction on $t$ gives

 \begin{eqnarray*}
  \E(e^{h(Z_1+ \cdots +Z_t)}) &=  &
  \E[\E(e^{h(Z_1+ \cdots+ Z_{t})} |e_1, \dots, e_{t-1} )]\\
&= &
\E[e^{h(Z_1+ \cdots+ Z_{t-1})}\E(e^{hZ_t}|e_1, \dots, e_{t-1} )]\\
&\leq &
\E[e^{h(Z_1+ \cdots+ Z_{t-1})}e^{h^2\eps \gamma_t}]\\
&\leq &
e^{h^2 \eps \sum_{i=1}^t \gamma_i}.
\end{eqnarray*}

Inserting this in  (\ref{eqn:concentration0}) we have

$$\Pr(X_t \ge n) \le e^{h^2 \eps \sum_{i=1}^t \gamma_i -hn}. $$

Since $\sum_{i=1}^t \gamma_i= O(n \log n)$ and $\eps=\log^{-1} n$, if
we set $h$ to be a sufficiently small positive constant, the right
hand side is $e^{-\Omega (n)} = n^{-\omega (1)}$, proving
(\ref{eqn:concentration}).
(Of course we could do much better---for
$\Pr(|X| > \lambda) <n^{-\omega(1)}$,
$\lambda > \Omega(\sqrt{n\log n })$
is enough---but this makes no difference for our purposes.)
That
$\Pr (X_n \le -n) < n^{-\omega (1)}$
is proved similarly, using $ -X_t$ in place of $X_t$.
(Actually we only use the bound on
$\Pr (X_n \geq n)$.)

\section{The properties $\B$ and $\r$}\label{BandR}

In this section we define the properties $B_i$ and $\r_i$
and observe that in dealing with (\ref{Ri}) and (\ref{Bi})
we can work with $G(n,p_i)$ rather than $G_i$,
where
$$p_i =
1- \frac{i}{{n \choose 2}}.$$

We will actually define properties $\B$ and $\r(p)$
($p\in [0,1]$); the events $\B_i$ and $\r_i$ are then
$$\{\mbox{$G_i$ satisfies $\B$}\}$$
and
$$\{\mbox{$G_i$ satisfies $\r(p_i)$}\}.$$

We need some notation.
For  a finite set $A$ and $\ww: A \rightarrow \Re^+$
($:=[0, \infty)$),
set
$$\overline \ww (A)= |A|^{-1} \sum_{a\in A} \ww(a),$$
$$\max \ww (A)= \max_{a \in A}\ww(a),$$
and
$$\maxr \ww(A)= \overline \ww (A)^{-1}\max \ww (A)$$
(the maximum normalized value of $\ww$), and write
$\med \ww(A)$ for the median of $\ww$.

%\begin{itemize}

%\item $\overline \ww (A)= \frac{1}{|A|} \sum_{a\in A} \ww(a)$ (mean
%value of $\ww$),

%\item  $\max \ww (A)= \max_{a \in A}$ (maximum value of $\ww$),

%\item $\maxr \ww(A)= \frac{ \max \ww (A)}{\overline \ww (A)}$
%(maximum normalized value of $\ww$),

%\item $\med \ww(A)$ is the median of $\ww$.

%\end{itemize}

For a graph $G$ on vertex set $V$
and $Z\in \C{V}{v}$, let $
%\ww(Z)=
\ww_G(Z) = \Phi(G-Z)$.  Thus, $\ww_{G}(Z)$ denotes the number of
$H$-factors in the graph induced by $G$ on the vertex set $V\backslash Z$.
We also use $\ww_G(K) = \ww_G(V(K))$
for $K\in\h(G)$.

The property $\B$ for a graph $G$ is
\begin{equation} \label{equation:B}
\B(G)=\{\maxr \ww_G (\h(G) ) =  O(1)\}.
\end{equation}

Thus, $\B(G)$ asserts that for any copy $K$ of $H$ in $G$,
the number of $H$-factors containing $K$ is not much more than the average.

There are two parts to
the definition of $\r(p)$.  The first refers to the following setup,
in which ($V=[n]$ and)
expectations refer to $G(n,p)$.

Given $A\sub V(H)$, $E'\sub E(H)\sm E(H[A])$,
$\psi$ an injection from $A$ to $ V$, and $G\sub K_n$,
let $X(G)$ be the number of injections
$\varphi:V(H)\ra V$ with
\begin{equation}\label{phipsi}
\mbox{$\varphi \equiv \psi$ on $A$}
\end{equation}
and
\begin{equation}\label{xyphi}
xy\in E' \Ra \varphi(x)\varphi(y)\in E(G).
\end{equation}
Write $X(G(n,p))$ in the obvious way as a polynomial in
variables $\ttt_e=\1_{\{e\in E(G(n,p))\}}$, $e\in E(K_n)$:
\begin{equation}\label{Xh}
X(G(n,p)) =h(\ttt) =\sum_{\varphi}\ttt_{\varphi(E')},
\end{equation}
where $\ttt =(\ttt_e:e\in E(K_n))$ and
the sum is over injections $\varphi:V(H)\ra V$
satisfying (\ref{phipsi}).
Then $h$ is multilinear, $O(1)$-normal and homogeneous of degree
$d=|E'|$, and we set
\begin{equation}\label{hE*}
\E^* =\max\{\E_Lh:|L|<d\}.
\end{equation}
(Note this includes $L=\0$, corresponding to $\E h$.)

For the second part of the definition, we use $D(p)$ for
the expected number of
copies of $H$
in $G(n,p)$ using a given $x\in V$; that is,
$$D(p)= v(n-1)_{v-1}p^m =\Theta(n^{v-1}p^m). $$

\begin{definition}\label{Rdef}
We say $G\sub K_n$
satisfies $\r(p)$ if the following two
properties hold.

{\rm (a)} For $A$, $E'$ and $\psi$ (and associated notation)
as above:
if $\E^* = n^{-\Omega(1)}$, then for any $\gb(n)=\go(1)$,
$X(G)< \gb(n)$ for large enough $n$;
if $\E^* \ge n^{-o(1)}$, then
for any fixed $\eps>0$ and large enough $n$,
$X(G) < n^{\eps}\E^*$.

{\rm (b)}  For each $x\in V$,
$|D(x,G) -D(p)|= o(D(p))$
\end{definition}

Recall $D(x,G)$ is the number of copies of $H$ (in $G$)
containing $x$. Thus (b) says that for any $x$,
the number of copies of $H$ containing $x$ is
close to its expectation.

We pause for the promised

\begin{lemma}  \label{proofofxibd}

For $i\leq T$ ($T$ as in Section \ref{Outline}),
$\B_{i-1}$ and $\r_{i-1}$
(that is, $\B$ and $\r(p_{i-1})$ for $G_{i-1}$) imply
{\rm (\ref{xibd})}.

\end{lemma}
{\em Proof.}
Write $\ww$ for $\ww_{G_{i-1}}$.
We first observe that
$\B$ and $\r_{i-1}$ imply that, for any $K \in \h(G_{i-1})$,
\begin{equation}\label{fracK}
\ww(K)/\Phi(G_{i-1}) = O(1/D(p_{i-1}))
\end{equation}
(the left side is the fraction of $H$-factors in $G_{i-1}$
that use $K$), since
\begin{eqnarray*}
\Phi(G_{i-1}) &=& \frac{v}{n}~\ww(\h(G_{i-1}))
%~\sum\{\ww(K'):K'\in \h(G_{i-1})\}
\\
&=& \frac{v}{n}~\Om(|\h(G_{i-1})|\max\ww(\h(G_{i-1})))
%\label{PG1}
\\
&=& \Om(D(p_{i-1})\ww (K))
%\label{PG2}
\end{eqnarray*}
(using $\B$ in the second line and part (b) of $\r_{i-1}$ in
the third).
%(using $\B$ for (\ref{PG1}) and part (b) of $\r_{i-1}$ for
%(\ref{PG2})).

On the other hand, a simple application of part (a) of $\r_{i-1}$
shows that, for any $e\in E(G_{i-1})$,
the number of $K\in \h(G_{i-1})$
containing $e$ is at most $\gb$ for some $\gb=\gb(n)$ satisfying
$\gb^{-1}D(p_{i-1}) =\go(\log n)$.
(Here we need the observation that for $i\leq T$,
$D(p_{i-1}) =\omega(\log n)$.)
Combining this with (\ref{fracK}) we have (\ref{xibd}).\qed

\bigskip
We also need the ``$p$-version" of $\A_t$:
\begin{equation}\label{A'}
\A(p) =\{\log |\f(G)| >\log|\f_0| - \sum_{i=1}^t\gc_i -O(n)\},
\end{equation}
where $t=\lceil (1-p)\C{n}{2}\rceil$
(and $\gc_i$ is as in (\ref{gci})).

According to Lemma~\ref{lemma:twomodels},
(\ref{Ri}) and (\ref{Bi})---and thus, as noted at the end of
Section~\ref{Outline}, Theorem~\ref{theorem:3'}---will
follow from the next two lemmas.
% in which we again
%use $\vartheta(n) = n^{-1/d(H)}(\log n)^{1/m}$.
\begin{lemma}\label{R'lemma}
For $p> \go(n^{-1/d(H)}(\log n)^{1/m})$,
$$\Pr(\mbox{$G(n,p)$ satisfies $\r(p)$}) = 1-  n^{-\go(1)}.$$
\end{lemma}
(The assumption on $p$ is only needed for (b) of Definition \ref{Rdef}.)
\begin{lemma}\label{B'lemma}
For $p> \go(n^{-1/d(H)}(\log n)^{1/m})$,
$$\Pr(\mbox{$G(n,p)$ satisfies $\A(p)\r(p)\bar{\B}$}) = n^{-\go(1)}.$$
\end{lemma}
These are proved in Section~\ref{Reg}
and Sections~\ref{PLB'} and \ref{PLC} respectively.

\section{Regularity}\label{Reg}

Here we verify Lemma~\ref{R'lemma}; that is,
we show that w.v.h.p. $G=G(n,p)$ satisfies (a) and (b)
of Definition \ref{Rdef}.

For (a),
since there are only
$n^{O(1)}$ possibilities for $A,E',\psi$, it's enough to
show that for any one of these the probability of violating
(a) is $n^{-\omega(1)}$.
This is immediate from the results of Section
\ref{Conc}:
recalling that each $(A,E',\psi)$
corresponds to a homogeneous, multilinear, $O(1)$-normal
polynomial $h$ as in (\ref{Xh}) and $E^*$ given by (\ref{hE*}),
we find that
the first part of (a) is given by Theorem
\ref{lastVu},
and the second by
Corollary \ref{Vucor} with $A = \frac{1}{2}n^{\eps}\E^*$
%%%%%%%%%%%%%%%%%%%%%%%%%%%%%%%%%%%%%%%%%%%%%
(and any $\eps$ smaller than the present one).
%%%%%%%%%%%%%%%%%%%%%%%%%%%%%%%%%%%%%%%%%%%%%

For (b), write $D(x,G)$ in the natural way
as a polynomial of degree $m$
(recall $m=|E(H)|$) in the variables
$t_e = \1_{\{e\in E(G)\}}$ ($e\in E(K_n)$);
namely (with $t=(t_e:e\in E(K_n))$),
$$
D(x,G) = f(t) :=\sum\{t_K:K\in \h_0(x)\},
$$
where $\h_0(x)$ is the set of copies of $H$ in $K_n$ containing $x$
(and, as in Section 5, $t_K=\prod_{e\in K}t_e$).

We have
\begin{equation}\label{appl1}
\E f= \Theta (n^{v-1} p^{m})=\go(\log n),
\end{equation}
and intend to show concentration of $f$ about its mean
using Theorem \ref{theorem:combination};
thus we need to say something about the expectations $\E_Lf$.

For $L\sub E(K_n)$ with $1\leq |L|=l <m$ we have
$$
\E_Lf = p^{m-l}N(L),
$$
where
$N(L)$ is the number of $K \in \h_0(x)$ with $L\sub E(K)$.
Let $W=V(L)\cup\{x\}$,
with $V(L)\sub V$
the set of vertices incident with edges of $L$,
and $v'=|W|$.
Then $N(L) =\Theta(n^{v-v'})$ if the graph
$H':= (W,L)$ is (isomorphic to) a subgraph of $H$,
and zero otherwise.

Thus, in view of (\ref{appl1}), we have
\begin{eqnarray*}
\E f/\E_Lf &= &\Om (n^{v'-1}p^l) =\Om (n^{[(v'-1)/l-(v-1)/m]l})\\
&=& \Om(n^{[1/d(H')-1/d(H)]l})
= n^{\Om(1)}
\end{eqnarray*}
(where we used strict balance of $H$ for the final equality).

Combining this with (\ref{appl1}), we have the hypotheses of
Theorem \ref{theorem:combination}, so also its conclusion,
which is (b).

\section{Proof of Lemma~\ref{B'lemma}}\label{PLB'}

We now assume $p$ is as in Lemma~\ref{B'lemma}
and slightly simplify our notation,
using $G$, $\r$ and $\A$ for $G(n,p)$, $\r(p)$
and $\A(p)$.
Events and the notation ``$\Pr$" now refer to $G$; so
for instance $\Pr(\A)$ is the probability that $G$ satisfies $\A$.

We will
get at $\B$ via an auxiliary event $\cee$.
Write $\h_0$ for the collection of $v$-subsets of $V$, and,
for $Y\sub V$ with $|Y|\leq v$,
$\h_0(Y)$ for $\{Z\in \h_0:Z\supseteq Y\}$.
We extend the weight function $\ww=\ww_G$ to such $Y$ by setting
$$
\ww(Y)= \sum \{\ww(Z): Z\in \h_0(Y)\}.
$$
Thus $\ww(Y)$ is the number of ``partial $H$-factors" of size
$\frac{n}{v}-1$ in $G-Y$.
The property $\cee$ for $G$ is
\begin{equation}\label{C}
%Y\in \C{V}{v-1} \Ra
\mbox{for any $Y\in \C{V}{v-1} $,
$\max \ww (\h_0(Y) )
\le \max \{n^{-2(v-1)} \Phi(G), 2 \med \ww ( \h_0(Y)) \}$}
\end{equation}

\noindent (where, recall,
$\Phi(G)$ is the number of $H$-factors in $G $).

Lemma~\ref{B'lemma} follows from the next two lemmas.
\begin{lemma}\label{Clemma}
$~~~\Pr(\A\r\bar{\cee}) = n^{-\go(1)}$
%$~~~\Pr(\mbox{$G$ satisfies $\A\r\bar{\cee}$}) < n^{-C}$.
\end{lemma}

\begin{lemma}\label{BClemma}
$~~~\Pr(\r\cee\bar{\B}) = n^{-\go(1)}$
%$~~~\Pr(\mbox{$G$ satisfies $\r\cee\bar{\B}$}) < n^{-C}.$
\end{lemma}

The proof of
Lemma~\ref{Clemma}, which is really the heart of
the matter,
is deferred to the next section.
Here we deal with Lemma~\ref{BClemma}.  We first need
to say that $\cee$ implies that
$
\maxr \ww(\h_0) = O(1),
$
or, equivalently, that for some positive constant $\gd$,
\begin{equation}\label{KHO}
|\{K\in \h_0:\ww(K) \geq \gd \max\ww(\h_0)\}| = \Omega(|\h_0|)
~(= \Omega(n^v)).
\end{equation}
This is an instance of a simple and general
deterministic statement that we may formulate
(more precisely than necessary) as follows.
Let $V$ be any set of size $n$
and $\ww:\C{V}{v}\ra \Re^+$.
%(For this lemma, the function $\ww$ is
% arbitrary and has nothing to do with the number of factors.)
For $X\sub V$ of size at most $v$ write $\h_0(X)$
for the collection of $v$-subsets of $V$ containing $X$,
and for simplicity set $\psi(X) =\max\ww(\h_0(X))$.
Let $B$ be a positive number.

\begin{lemma}\label{deterministic}
Suppose that for each $Y\sub V$ satisfying $|Y|=v-1$ and
$\psi(Y)\geq B$ we have
$$
|\{Z\in \h_0(Y):\ww(Z)\geq \frac{1}{2} \psi(Y)\}|\geq \frac{n-v}{2}.
$$
Then for any $X\sub V$ with $|X|=v-i$ and $\psi(X)\geq  2^{i-1}B$
we have
\begin{equation}\label{Ni}
|\{Z\in \h_0(X):\ww(Z)\geq \frac{1}{2^{i}} \psi(X)\}|\geq
\left(\frac{n-v}{2}\right)^i\frac{1}{(i-1)!}~.
\end{equation}
\end{lemma}
{\em Proof.}
Write $N_i$ for the r.h.s. of (\ref{Ni}).
We proceed by induction on $i$, with the case $i=1$ given.
Assume $X$ is as in the statement and choose $Z\in \h_0(X)$
with $\ww(Z)$ maximum (i.e. $\ww(Z)=\psi(X)$).
Let $y\in Z\sm X$ and $Y=X\cup\{y\}$.
Then $|Y|=v-(i-1)$ and
$\psi(Y) =\psi(X)\geq 2^{i-1}B$ ($\geq 2^{i-2}B$);
so by our induction hypothesis
there are at least $N_{i-1}$ sets
$Z'\in \h_0(Y)$ with
$\ww(Z')\geq 2^{-(i-1)}\psi(Y)$ ($=2^{-(i-1)}\psi(X)$).
For each such $Z'$, $Z'\sm\{y\}$ is a $(v-1)$-subset of $V$ with
$\psi(Z'\sm\{y\})\geq \ww(Z')\geq B$.
So (again, for each such $Z'$) there are at least $(n-v)/2$
sets $Z''\in\h_0(Z'\sm\{y\})$ with
$$
\ww(Z'')\geq \psi(Z'\sm\{y\})/2\geq 2^{-i}\psi(X).
$$
The number of these pairs $(Z',Z'')$ is thus at least
$N_{i-1}(n-v)/2$.
On the other hand, each $Z'$ associated with a given $Z''$ is $Z''\sm\{u\}\cup \{y\}$ for some $u\in Z''\sm (X\cup\{y\})$;
so the number of such $Z'$ is at most $i-1$ and the lemma
follows.\qed

{\em Proof of Lemma} \ref{BClemma}.
%%%%%%%%%%%%%%Added +1 to the exponent to make up for $n-v$ vs $n$
Set $\gc = [2^{v+1} (v-1)!]^{-1}$ and $\gd = 2^{-v}$.
Now $\cee$ implies the hypothesis
of Lemma \ref{deterministic} with
$B= (2n)^{-(v-1)}\Phi(G)$
%%%%%%%%%%%%%%%%%%%%%%%%%%%%%%%%%%%%%%%%%%%%%
(actually with any $B $ greater than $ n^{-2(v-1)}$),
%%%%%%%%%%%%%%%%%%%%%%%%%%%%%%%%%%%%%%%%%%%%%
and we have (trivially)
$$\psi(\0) \geq n^{-(v-1)}\Phi(G) = 2^{v-1}B.$$
Thus (\ref{Ni}) applies, yielding (\ref{KHO}), now
with constants specified:
$$
|\{K\in \h_0:\ww(K)\geq\gd \max\ww(\h_0)\}| > \gc n^v.
$$

Let $\zz$ be the largest power of 2 not exceeding $\max\ww$
and
$$\Z= \{Z\in \h_0:\ww(Z)>\gd \zz\}.$$
For $X\sub V$ with $|X|\leq v$ let
$\Z(X)= \{Z\in\Z:X\sub Z\}$,
and say $X\sub V$ with $|X|\leq v$ is
{\em good} if
$|\Z(X)| >\gc n^{v-|X|}$.
Thus in particular,
\begin{equation}\label{Cimplies}
\mbox{{\em the empty set is good}}
\end{equation}
(as are all singletons, but we don't need this).

Fix an ordering $a_1\dots a_v$ of $V(H)$.
For distinct $x_1\dots x_r\in V$,
write $S(x_1\dots x_r)$ for the collection of copies $\varphi$
of $H$ in $K_n$ for which
$$\gp (a_i)=x_i ~~~\mbox{for $i\in [r]$,}$$
$$\gp(a_i)\gp(a_j)\in E(G)
~~\mbox{whenever $i,j\geq r$ and $a_ia_j\in E(H)$,}$$

and $\gp(V(H))\in\Z$.

For $r\in \{0,\dots v\}$ let
$N_r =N(a_r)\cap \{a_{r+1}\dots a_v\}$ and $d_r = |N_r|$.
(In particular $d_v=0$.)
Let $\Y(x_1\dots x_r) $ be the event
$$\{|S(x_1\dots x_r)| =
\Omega(p^{d_r+\cdots +d_{v-1}}n^{v-r})\}.$$
%$$\Y(x_1\dots x_r) =
%\{|\Z(x_1\dots x_r)| =
%\Omega(p^{d_r}\sum_{x_{r+1}}|\Z(x_1\dots x_{r+1})).$$
In particular
$$S(\0) = \{\varphi\in \h(G): \ww(\varphi(V(H)))>\gd\zz\}$$
and
$\Y(\0)$  is the event
$$\{|S(\0)| = \Omega(p^mn^v)\}.$$

For
(distinct) $x_1\dots x_r\in V$ let $Q(x_{1} \dots x_{r})$  be the event
$$\{\mbox{$\{x_1\dots x_r\}$ is good}\}
\wedge \ov{\Y}(x_1\dots x_r).$$

Since $\bar{\B}\cee\sub Q(\0)$ (by (\ref{Cimplies})),
Lemma~\ref{BClemma} will follow if we can show
$$\Pr (\r Q(\0)) = n^{-\go(1)}. $$
For inductive purposes we will actually prove the
more general statement that
for any $r $ and $x_1\dots x_r$,
\begin{equation}\label{ifgood}
\Pr(\r Q(x_1\dots x_r))= n^{-\go(1)}.
\end{equation}
{\em Proof of} (\ref{ifgood}).
We proceed by induction on $v-r$.
The case $r=v$ being trivial
(since for $v$-subsets of $V$ being good is the same as belonging
to $\Z$),
we consider $r<v$, letting $X=\{x_1\dots x_r\}$.

Let $\pee$ be the event
$$\{y\in V\sm X, \mbox{$X\cup \{y\}$ good}
~\Rightarrow ~\Y(x_1\dots x_r,y)\}.$$
By inductive hypothesis it's enough to show
\begin{equation}\label{rcpQ}
\Pr(\r\pee Q(x_1\dots x_r)) < n^{-\go(1)}
\end{equation}
(since
$\Pr (\r Q) \le \Pr (\r \pee Q) + \Pr (\r \ov \pee) $
and, by induction,
$\Pr (\r \ov \pee) = n^{-\go(1)} $).

Note that if $X$ is good then
\begin{equation}\label{manyx}
|\{y:\mbox{$X\cup\{ y\}$ good}\}| = \Omega(n).
\end{equation}
We need to slightly relax $\r$ to arrange that
the edges between $x_r$ and $V\sm X$ are
independent of our conditioning.
Say $G$ satisfies $\r_X$ if it satisfies (a) in the
definition of $\r$ ($=\r(p)$) whenever $A=\{a_1\dots a_r\}$,
$\psi(a_i)=x_i$ ($i\in [r]$) and $E'\sub E(H-A)$.

If
$\r\pee \wedge \{\mbox{$X$ good}\}$ holds
but $\Y(x_1\dots x_r)$ does not, then we
have the following situation.
There is some $J=2^k$ with $k$ an integer not exceeding $n\log n$
(see (\ref{mg0})) so that (with $\Z$, ``good" and other quantities
as in the preceding discussion, but now defined in terms of
this $J$)

\vskip2mm

\hskip10mm
(i) $\r_X$ holds;

\hskip10mm
(ii) there are at least $\Omega( n)$ $y$'s in $V\sm X$ for which
we have $\Y(x_1\dots x_r,y)$
 (by (\ref{manyx})); but

\hskip10mm
(iii) $\Y(x_1\dots x_r)$ does not hold.

\vskip2mm
Note that, for a given $J$, properties (i) and (ii) depend only on
$G':=G-X$.
Since the number of possibilities for $J$ is at most $n\log n$,
it is thus enough to show that for any $J$
and $G'$ satisfying (i) and (ii)
(with respect to $J$),
\begin{equation}\label{showX}
\Pr(\ov{\Y}(x_1\dots x_r)|G') = n^{-\go(1)}.
\end{equation}

Given such a $G'$
(actually, any $G'$),
$|S(x_1\dots x_r)|$ is naturally expressed as
a multilinear polynomial in the variables
$$t_u:= \1_{\{x_ru\in E(G)\}}  \,\,\,\,\,\,\, u\in V\sm X;$$
namely
$$|S(x_1\dots x_r)| = g(t):=\sum_U\ga_{_U}t_{_U},$$
where $U$ ranges over $d_r$-subsets of $V\sm X$, and
$\ga_{_U}$ is the number of copies $\psi$ of
$K:=H-\{a_1\dots a_r\}$ in $G'$ with
$$\psi(N_r)= U$$
and
\begin{equation}\label{psiZ}
\psi(\{a_{r+1}\dots a_v\})\cup X\in \Z.
\end{equation}
In order to apply Theorem~\ref{theorem:combination}
we normalize and consider
$$f(t) = \ga^{-1} g(t),$$
with $\ga $ the maximum of the $\ga_{_U}$'s.
We then need the hypotheses of Theorem~\ref{theorem:combination};
these will follow (eventually) from $\r_X$.

We may rewrite
$$
g(t) = \sum_{y\in V\sm X}
\sum\{t_{\varphi(N_r)}:\varphi\in S(x_1\dots x_r,y)\}.
$$
Noting again that the variables $t_u$ ($u\in V\sm X$) are
independent of $G'$ (which determines the sets
$S(x_1\dots x_r,y)$), and using property (ii),
we have
\begin{equation}\label{Eg}
\E g =  p^{d_r} \sum_{y\in V\sm X} |S(x_1\dots x_r,y)|
=\Omega(p^{d_r +\cdots +d_{v-1}}n^{v-r}).
\end{equation}
Note that if $d_r=0$ then there is nothing random
at this stage and we have
$$
|S(x_1\dots x_r)| = \sum_{y\in V\sm X} |S(x_1\dots x_r,y)|=
\Omega(p^{d_r +\cdots +d_{v-1}}n^{v-r}),
$$
which is what we want; so we may assume from now on that $d_r>0$.

If we set $H' = H-\{a_1\dots a_{r-1}\}$ then
$d_r+\cdots +d_{v-1}= e(H')$ and $v-r=v(H')-1$, so that
the r.h.s. of (\ref{Eg}) is
$\Omega(p^{e(H')}n^{v(H')-1})$.
Thus, since $H$ is strictly balanced,
we have
\begin{equation}\label{Ef'''}
\E g = \left\{\begin{array}{ll}
\omega (\log n)  & \mbox{if $r=1$}\\
n^{\Omega(1)}&\mbox{if $r>1$.}
\end{array}\right.
\end{equation}

We are going to prove that
\begin{equation}\label{Ef}
\E f = \omega (\log n)
\end{equation}
(in most cases it will be $n^{\Omega(1)}$)
and
\begin{equation}\label{Ef2}
\max\{\E_T f: T\sub V\sm X, 0<|T|<d_r\} = n^{-\Omega(1)}\E f.
\end{equation}

If we have these then Theorem~\ref{theorem:combination}
says that w.v.h.p. $f$ and (therefore) $g$ are close to
their expectations, which in view of (\ref{Eg}) is what we want.

For the proofs of (\ref{Ef}) and (\ref{Ef2}),
it will be more convenient to work with the partial derivatives of $g$ than with those of $f$.
As in Section 5
we use $ t_e=\1_{\{e\in E(G)\}}$,
$t_S=\prod_{e\in S}t_e$ and
\begin{equation}\label{sss}
\ttt = (t_e:e\in E(K_{V\sm X}))
\end{equation}
(where, recall, $K_{V \sm X}$ is the complete graph on $V \sm X$).

Let $T\sub V\sm X$ with $0<\lll:=|T|\leq d_r$.
(Note we now include $\lll=d_r$, in which case $\E_Tg$ is just $\ga_T$.)

Since we are only interested in upper bounds on the partial
derivatives of $g$, we may now disregard the requirement
(\ref{psiZ}); thus we use
\begin{equation}\label{h(s)}
p^{-(d_r-\lll)}\E_T g
\leq h(\ttt) := \sum_{\varphi}t_{\varphi(E(K))},
\end{equation}
where the sum is over injections
\begin{equation}\label{phis}
\mbox{$\varphi:V(K)\ra V\sm X$
with $\varphi(N_r) \supseteq T$}
\end{equation}
(and, again, with the obvious meaning for $\varphi(E(K))$).

Set $\E^* =\max\{\E_Lh:L\sub E(K_{V\sm X}),|L|<|E(K)|\}$.
We assert that there is a positive constant $\eps$
(depending only on $H$) so that (for large enough $n$),
\begin{equation}\label{E*}
p^{d_r-\lll}\E^*  < n^{-\eps}\E g.
\end{equation}
%(Of course we have already randomized $G'$,
%but we now need to reexamine this randomization
%to see what $\r_X$ says about the $\E_Tg$'s.)
%%%%%%%%%%%%%%%%%%%%%%%%%%%%%%%%%%%%%%%%%%%%%%%%%%%%%
(Of course we have already chosen $G'$,
but we now need to reexamine the randomization
that produced it
to see what $\r_X$ says about the $\E_Tg$'s.)
%%%%%%%%%%%%%%%%%%%%%%%%%%%%%%%%%%%%%%%%%%%%%%%%%%%%%

Before proving (\ref{E*}), we show that it gives
(\ref{Ef}) and (\ref{Ef2}).
For (\ref{Ef}) we need
$$
\ga^{-1}\E g= \omega(\log n).
$$
We apply(\ref{E*}) with $T=U$, a $d_r$-subset of $V\sm X$
(in which case,
as noted above, $\ga_U$ is just $\E_Ug$).
We consider two possibilities.
If $\E g \geq n^{\eps/2}$ then $\r_X$ gives
(note here $d_r-\lll=0$)
$$\ga_U < n^{\eps/4}\max\{1,\E^*\} \le  n^{-\eps/4}\E g.
$$
Otherwise we have $\E^* < n^{-\eps/2}$.  In this case,
recalling that $\E g = \omega(\log n)$ (see (\ref{Ef'''})),
we can choose $\gb(n)=\go(1)$ with
$\gb(n)^{-1}\E g = \go(\log n)$,
and $\r_X$ guarantees that $\ga_U\le \gb(n)$.
So in either case we have $\ga^{-1}\E g = \go(\log n)$.
Since none of these bounds depended on the choice
of $U$, this gives (\ref{Ef}).

For (\ref{Ef2}) we need $\E_Tg = n^{-\Omega(1)}\E g$
for any $T$ as in (\ref{Ef2}).
This follows from $\r_X$ if
we assume, as we may, that the constant $\eps$ in (\ref{E*}) is
less than $ 1/d(H)$, since we then have, e.g.,
$$
E_Tg <
p^{d_r-\lll}n^{\eps/2}\max\{1,\E^*\} \le n^{-\eps/2}\E g.
$$

For the proof of (\ref{E*}),
fix $L\sub E(K_{V\sm X})$ and
let $k = |E(K)|-|L|$.
We have $\E_Lh = p^kN_L$, where
$N_L$ is the number of $\varphi$ as in
(\ref{phis}) with $\varphi(E(K))\supseteq L$.
In particular each such $\varphi$ satisfies
$\varphi(V(K))\supseteq W:=T\cup V(L)$,
where, as earlier,
$V(L)\sub V\sm X$ is the set of vertices incident with edges of $L$.
Letting $W =\{w_1\dots w_s\}$, we have
$N_L= \sum N_L(b_1\dots b_s)$,
where $(b_1\dots b_s)$ ranges over $s$-tuples of distinct
elements of $V(K)$ and
the summand is the number of $\varphi$'s as above
with $\varphi(b_i) =w_i$ for $i\in [s]$.
Since there are only $O(1)$ choices for the $b_i$'s, we will have
(\ref{E*}) if we show (for any choice of $b_i$'s)
\begin{equation}\label{E**}
p^{d_r-\lll+k} N_L(b_1\dots b_s) = n^{-\Omega(1)} \E g.
\end{equation}
Let (given $b_i$'s)
$H'' =H[\{a_r,b_1\dots b_s\}]$.  Then
$$N_L(b_1\dots b_s) < n^{v-r-s} = n^{v-r-(v(H'')-1)},$$
while $k\geq \lll+d_{r+1}+\cdots + d_{v-1} -e(H'')$
(since $|E(K)|=d_{r+1}+\cdots +d_{v-1}$ and
$E(H'')$ contains $\varphi^{-1}(L)$ and at least $\lll$ edges
joining $a_r$ to $V(K)$).
Thus
$$p^{d_r-\lll+k} N_L(b_1\dots b_s) <
p^{d_r+\cdots +d_{v-1}}n^{v-r} [n^{v(H'')-1}p^{e(H'')}]^{-1}.$$

Since $H$ is strictly balanced (and since $v(H'')-1 \geq \lll> 0$),
the expression in square brackets is $n^{\Omega(1)}$,
so that, in view of (\ref{Eg}), we have (\ref{E**}).

This completes the proof of Lemma~\ref{BClemma}.

\section{Proof of Lemma~\ref{Clemma}}\label{PLC}

We continue to use $G$, $\r$ and $\A$ as in the preceding section.

Assume that we have $\A$ and $\r$ and that $\cee$ fails at $Y$
(see (\ref{C})).
Let $R= Y\cup \{x\}\in \h_0(Y)$ satisfy
$\ww(R) = \max \ww (\h_0(Y))$.
Note that
\begin{equation}\label{Kx}
\ww(R) > n^{-2(v-1) } \Phi(G).
\end{equation}

Choose
$y\in V\sm Y$ with
$\ww(Y\cup \{y\}) \leq  \med \ww (\h_0(Y) )$, and
with $h(y,G-R)$ maximum subject to this restriction
(see Section \ref{Ent} for $h$),
and set $S=Y\cup \{y\}$.
Thus in particular
$$\ww(R) > 2\ww(S).$$

Note that the lower bound in (\ref{A'}) is now
$$
\frac{v-1}{v}n \log n + \frac{mn}{v} \log p -O(n)
$$
(see (\ref{mg0}) and (\ref{Em})),
while $D=D(p)$ (as in Definition \ref{Rdef}) satisfies
$$
\log D =(v-1)\log n + m\log p -O(1).
$$
Thus $\A$ and (\ref{Kx}) give
$$
\log \Phi((G-R) \geq \frac{v-1}{v}n \log n + \frac{mn}{v} \log p -O(n) ,
$$
while from $\r$
we have, using Lemma~\ref{Shearer},
$$
\log \Phi(G-R) \leq \frac{n}{2v}[\log ((1+o(1)) D) +h(y,G-R)].
$$
%%%%%%%%%%%%%%%%%%%%%%%%%%%%%%%%%%%%%%%%%%%%%%%%%%%%%%
(In more detail:
Lemma~\ref{Shearer} gives
$\log \Phi(G-R) \leq \frac{n-v}{v}\sum_{z\in V\sm R}h(z,G-R)$;
and we have $h(z,G-R)\leq h(y,G-R)$ for at least half
the $z$'s in $V\sm R$,
and $h(z,G-R) \leq \log D(z,G-R) < \log((1+o(1))D)$ for {\em every} $z$,
using $\r$ and the fact that for any random variable $\kappa$,
$H(\kappa)\leq \log |{\rm range}(\kappa)|$.)
%%%%%%%%%%%%%%%%%%%%%%%%%%%%%%%%%%%%%%%%%%%%%%%%%%%%%%%

Combining (and rearranging) we have (again using $\r$)
\begin{equation}\label{handD}
h(y,G-R) > (v-1)\log n + m\log p -O(1) > \log D_{G-R}(y) -O(1).
\end{equation}

Let $W=V\sm (Y\cup\{x,y\})$ and for $Z\in \C{W}{v-1}$ set
$$\ww'(Z) = \Phi(G-(W\cup Z)).$$
The corresponding weight functions $\ww_y$ on $\h(y,G-R)$
and $\ww_x$ on $\h(x,G-S)$
are
$$
\mbox{$\ww_y(K)= \ww'(V(K)\sm\{y\})~~$
and $~~\ww_x(K)= \ww'(V(K)\sm\{x\})$.}
$$
Thus $X(y,G-R)$ (see Section \ref{Ent}) is chosen
according to the weights $\ww_y$, and similarly
for $X(x,G-S)$.
Note also that $\ww_y(\h(y,G-R)) = \ww(R)$
and $\ww_x(\h(x,G-S)) = \ww(S)$.

According to Lemma~\ref{entlemma}, (\ref{handD})
implies that there are
$a,b\in {\rm range}(\ww_y)$ ($={\rm range}(\ww')$)
as in (\ref{ab}) for which
$J := \ww_y^{-1}([a,b])
%\cap \h(y,G-R)
$ satisfies
\begin{equation}\label{mostof}
|J| > \Omega(|\h(y,G-R)|)
\end{equation}
and
%\begin{equation}\label{mostofw}
$$
\ww_y(J) > .7\ww_y(\h(y,G-R))=.7\ww(R).
$$
%\end{equation}
Setting $J'=\ww_x^{-1}([a,b])
%\cap \h(x,G-R)
$,
we thus have
\begin{equation}\label{wyJ}
\ww_y(J)>.7\ww(R) ,
\end{equation}
while
\begin{equation}\label{wyJ'}
\ww_x(J')\leq \ww(S) < .5\ww(R).
\end{equation}

Once we condition on the value of $G[W]$,
$\ww_y(J)$ and $ \ww_x(J')$ are naturally expressed
as evaluations of a multilinear polynomial in variables
$\{t_u:u\in W\}$, as follows.
Given $U\sub W$ with $|U|\leq v-1$,
let $G^*_U$ be the graph obtained
from $G[W]$ by adjoining a vertex $w^*$ with neighborhood $U$.
Let $\k_U$ be the set of copies of $H$ in $G^*_U$
containing $\{w^*u:u\in U\}$, and
\begin{equation}\label{alphaU}
\ga_{_U} = \sum\{\ww'(V(K)\sm\{w^*\}):
K\in \k_U,\ww'(V(K)\sm\{w^*\})\in[a,b]\}.
\end{equation}
The desired polynomial is then
$$
g(t) =\sum_{U\sub W}\ga_{_U}t_{_U},
$$
and $\ww_y(J)$ and $\ww_x(J')$ are $g$ evaluated at
$\ttt':=\1_{\{z\in W:yz\in E(G)\}}$ and
$\ttt^{''}:=\1_{\{z\in W:xz\in E(G)\}}$.

Now $\r$
implies $|\h(y,G-R)| =\Theta(n^{v-1}p^m) =\go (\log n)$.
(In more detail:  $\r$ implies that $D_G(y) =\Theta(n^{v-1}p^m)$
and, as is easily seen, that the number of copies of $H$
in $G$ containing $y$ and meeting $R$ is $o(n^{v-1}p^m)$.)
Thus by (\ref{mostof}) we also have
$$|J| =\Theta(n^{v-1}p^m),$$
and (since $b<O(a)$)
\begin{equation}\label{wyJ*}
\ww_y(J) =\Theta(bn^{v-1}p^m).
\end{equation}
The key idea is now that we can use Corollary~\ref{Vucor2}
to deduce that, unless something unlikely
has occurred,
$\E g$ must be of a similar size.

Fix $\RRR\sub W$, say with $|\RRR|=\lll<v$.
For $d=\lll\dots v-1$, and
$\ttt$ as in (\ref{sss}) with $V\sm X$ replaced by $W$,
we consider the polynomial
\begin{equation}\label{hd}
h_d(\ttt) = \sum_z\sum_{\varphi} \ttt_{\varphi(E(H-z))},
\end{equation}
where $z$ ranges over vertices of $H$ of degree $d$
and $\varphi$ over injections $V(H)\sm\{z\} \ra W$ with
$\varphi(N_z)\supseteq \RRR$.
Then
$$
\ga_{\RRR} \leq b\cdot h_{\lll}(\ttt)$$
(``$\leq$" because of the restriction
$\ww'(V(K)\sm\{w^*\})\in[a,b]$ in
(\ref{alphaU}))
and
$$
\E_{\RRR}'g\leq b \sum_{d>\lll}p^{d-\lll}h_d(\ttt).
$$

Let
$$
\E_d^* =\max\{\E_Lh_d:L\sub E(K_W),|L|< m-d\}.
$$

We assert that there is a positive constant $\eps$
(depending only on $H$) so that (for each $d$),
\begin{equation}\label{E*'}
p^{d-\lll}\E_d^*  < n^{-\eps}n^{v-1}p^m.
\end{equation}
The proof of this is essentially identical to that of
(\ref{E*}) and we omit it.
(Actually (\ref{E*'}) is contained in (\ref{E*}):
it's enough to prove (\ref{E*'}) with $\E_d^*$
replaced by $\E^*(z)$ gotten by replacing $h_d$ in
the definition of $\E^*_d$ by the inner sum in (\ref{hd});
but this inner sum is bounded by the polynomial $h(\ttt)$
in (\ref{h(s)}) with $r=1$,
$a_1=z$ and $d_r = d$
(``bounded by" rather than ``equal to" because we've replaced
$V\sm \{x_1\}$ by the slightly smaller $W$);
and, finally, note that we actually proved
(\ref{E*}) with $\E g$ replaced by
$\Omega(p^{d_r +\cdots +d_{v-1}}n^{v-r})$ (see (\ref{Eg}).)

Now intending to apply the results of Section \ref{Conc},
we consider
$f=\ga^{-1}g$, where $\ga=\max_{_U}\ga_{_U}$.
Then
$$f(\ttt') = \ga^{-1}\ww_y(J) =
\Theta(\ga^{-1} b n^{v-1}p^m)$$
(see (\ref{wyJ*})), and we have
\begin{equation}\label{Ef'}
 f(\ttt') =\omega (\log n)
\end{equation}
and
\begin{equation}\label{Ef2'}
\max\{\E'_T f: T\sub W, T\neq \0\} = n^{-\Omega(1)}f(\ttt').
\end{equation}
These are derived from (\ref{E*'}) in the same way as
(\ref{Ef}) and (\ref{Ef2}) were derived from (\ref{E*}),
and we will not repeat the arguments.

\medskip
In summary, $\A\r\bar{\cee}$ implies that there are $Y,x,y$
and (with  notation as above)
$a,b\in  {\rm range}(\ww')$
for which we have (\ref{Ef'}),
(\ref{Ef2'}), and, from (\ref{wyJ}) and (\ref{wyJ'}),
\begin{equation}\label{JJ'}
f(\ttt^{''}) < .8 f(\ttt').
\end{equation}

But for given $Y,x,y,a$ and $b$,
$f$ depends only on $G[W]$.  On the other hand,
given $G[W]$, $\ttt$ and $\ttt'$ are independent
r.v.'s, each with law ${\rm Bin}(W,p)$ (where
we say $\ttt=(\ttt_w:w\in W)$ has ``law ${\rm Bin}(W,p)$"
if the $\ttt_w$'s are independent, mean $p$ Bernoullis).
Thus the following simple consequence of
Theorem~\ref{inhomog} and Corollary~\ref{Vucor2} applies.
%\begin{claim}\label{finalclaim}

\mn
{\bf Claim.}
{\em
For any $\eps>0$ and d the following holds.
If $f$ is a
multilinear, normal polynomial of degree at most $d$
in n variables, $\gz(n) =\omega(\log n)$,
and $\ttt'$, $\ttt''$ are independent, each with law
${\rm Bin}([n],p)$, then
$$\Pr(f(\ttt') >
\max\{\gz(n),n^{\eps} \max_{T\neq\0}\E_T' f, (1+\eps)f(\ttt'')\})
=n^{-\omega(1)}.$$}
%\end{claim}
Since there are only polynomially many possibilities for
$Y,x,y,a$ and $b$, this gives Lemma~\ref{Clemma}.

\mn
{\em Proof of Claim.}
Set
$$A=\frac{1}{2}\max\{\gz(n),n^{\eps} \max_{T\neq\0}\E_T'f\}.$$
If $\E f\leq A$ then Corollary~\ref{Vucor2} gives
$$\Pr(f(\ttt')> A)= n^{-\omega(1)};$$
otherwise, by Theorem~\ref{inhomog},
\begin{eqnarray*}
\Pr(f(\ttt')>(1+\eps)f(\ttt'')) & <&
\Pr(\max\{|f(\ttt')-\E f|,|f(\ttt'')-\E f|\}> (\eps/3) \E f) \\
&=  & n^{-\omega(1)}.
\end{eqnarray*}

\section{Extensions} \label{extension}

As mentioned earlier, we will not repeat the above arguments
for general graphs or for hypergraphs; but we do want to stress
here that ``repeat" is the right word:
extending the proof of
Theorem \ref{theorem:3}
to produce (the counting versions of)
Theorems \ref{theorem:1-1}, \ref{theorem:hyper}
and \ref{theorem:1-1hyper} involves nothing but some minor
formal changes.  To elaborate slightly:

While strict balance is used repeatedly in the proof of
Theorem \ref{theorem:3}, a little thought shows that
all these uses are of the same type:
to allow us to say that, for some proper subgraph $H'$ of $H$,
and $p$ in the range under consideration,
\begin{equation} \label{essential}
n^{v(H')-1} p^{e(H')} = n^{\Omega (1)},
\end{equation}
as follows from
\begin{equation}\label{useofsb}
e(H')/(vH')-1) < m/(v-1)
\end{equation}
(this is what we get from strict balance) and the fact that $p\geq
n^{-1/d(H)}$ (of course it is actually somewhat bigger). See the
last few paragraphs of Section \ref{Conc}.

For general graphs we no longer have
(\ref{useofsb}) but can still guarantee
(\ref{essential})---though, of course,
at the cost of some precision in the results---by taking
$p \ge n^{-1/d^{\ast}(H) +\epsilon}$ for some fixed positive
$\epsilon$.  This leads to
Theorem \ref{theorem:1-1}.

The extensions to hypergraphs are similarly unchallenging, basically
amounting to the observation that our arguments really make no use
of the assumption that edges have size two.

\section{Appendix A: equivalence of Theorems \ref{theorem:2} and
\ref{theorem:3}}\label{Appendix}

Of course we only need to show
Theorem \ref{theorem:3} implies
Theorem \ref{theorem:2}.

Set
$\vartheta(n) =n^{-1/d(H)} (\log n)^{1/m}$
and $\mu(n,p) = (n^{v-1} p^m)^{n/v}$, and
write $G$ for $G(n,p)$ (for whatever $p$ we specify).

Given $K$,
we would like to show that there is a $C_K$ for which
$p=p(n)>C_K\vartheta(n) $ implies
\begin{equation}\label{Cworks}
\Pr(\Phi(G) \leq \mu(n,p)e^{-C_Kn})\leq n^{-K} ~~~\forall n.
\end{equation}
Suppose this is not true and for each $C,n$
define $g_C(n)$ by:
if $p=p(n)=g_C(n)\vartheta(n)$ then
$$\Pr(\Phi(G(n,p))\leq \mu(n,p)e^{-Cn}) = n^{-K}.$$

If the sequence $\{g_C(n)\}_n$ is bounded for some $C$
then we have (\ref{Cworks}) with $C_K$ the maximum of $C$ and
the bound on $g_C$.

So we may assume that $\{g_C(n)\}_n$ is unbounded for every $C$.
We can then choose $n_1<n_2< \cdots$ so that
$g_C(n_C)>C$ for $C=1,2,\ldots$.

Define the sequence $\{g(n)\}$ by
$g(n) =g_C(n_C)$, where $C$ is minimum with $n_C\geq n$.
Since $g(n)\ra \infty$,
Theorem \ref{theorem:3} says that there is a $C^*$ so that if
\begin{equation}\label{p(n)}
p=p(n)=g(n)\vartheta(n)
\end{equation}
then
$$\Pr(\Phi(G(n,p))\leq \mu(n,p)e^{-C^*n}) = n^{-\go(1)}.$$
Now choose $n_0$ so that for $n>n_0$ and $p$ as in
(\ref{p(n)})
$$\Pr(\Phi(G(n,p))\leq \mu(n,p)e^{-C^*n}) < n^{-K},$$
and then $C>C^*$ with $n_C>n_0$.
Then for $p$ still as in
(\ref{p(n)}) and $n=n_C$ we have
\begin{equation}\label{contradiction}
\Pr(\Phi(G(n,p))\leq \mu(n,p)e^{-Cn})<
\Pr(\Phi(G(n,p))\leq \mu(n,p)e^{-C^*n}) < n^{-K}.
\end{equation}
But this contradicts our definition of $g$, according to
which the l.h.s. of (\ref{contradiction}) is
$n^{-K}.$%\qed

\end{document}